# Non Parallelism of Ricci tensor of Pseudo-cylindric Metrics

A. Raouf Chouikha

**Abstract**

We descrive examples of metrics in the conformal class $[g]$ on complete conformally flat Riemannian manifolds $(M, g)$. These metrics have a constant scalar curvature and an harmonic curvature with non parallel Ricci tensor.[1]

## 1 Introduction

It is wellknown that the class $\mathcal{P}$ of compact Riemannian manifolds with parallel Ricci tensor is included in the class $\mathcal{C}$ of manifolds with constant scalar curvature.
Let $\mathcal{H}$ the class of such manifolds with harmonic curvature. The following inclusion holds
$$\mathcal{P} \subset \mathcal{H}.$$

Moreover, this inclusion is strict as Riemannian metrics given by A. Derdzinski [B], [D] and A. Gray [G].
In other words, there exist metrics with harmonic curvature and Ricci tensor non parallel. But examples of such metrics are very few.

Consider a Yamabe metric $g_0$ on a n-dimensional Riemannian manifold $(M, g)$. That means there is a $C^\infty$ positive function $u$ satisfying a differential equation such that $g_0 = u^{\frac{4}{n-2}} g$, so that its scalar curvature is constant. We are interested on the following problem:
When can $g_0$ have a harmonic curvature and non parallel Ricci tensor ?

---

[1] 2000 *Mathematics Subject Classification*, 53C21, 53C25, 58G30.

This concerns specially the conformally flat manifolds.

Let us consider the Riemannian cylindric product $(S^1 \times S^{n-1}, dt^2 + d\xi^2)$, where $S^1$ is the circle of length $T$ and $(S^{n-1}, d\xi^2)$ is the standard sphere. This metric is known to have a parallel Ricci curvature tensor . Furthermore, the number of Yamabe metrics is finite in the conformal class of the Riemannian product metric $[dt^2 + d\xi^2]$, [Ch], [M-P] .
We call a pseudo-cylindric metric any non trivial metric of the conformal class.
We also consider the manifold $M = S^n - \Lambda_k$, where $\Lambda_k$ is a finite point-set of the standard sphere $(S^n, d\xi^2)$.

We show that unless the trivial ones, all the pseudo-cylindric metrics in the conformal class $[dt^2 + d\xi^2]$ as well as the Yamabe metric of $(M, d\xi^2)$ belong $\mathcal{H}$ but not belong $\mathcal{P}$,

## 2  The Derdzinski metrics and the harmonic curvature

Let us consider a compact Riemannian manifold $(M, g)$, $\mathcal{R}$ is the curvature associated to the metric $g$, $r = Ric(g)$ is its Ricci tensor, and $D$ is the Riemannian connexion. We denote by $\delta\mathcal{R}$ its formal divergence. This curvature is harmonic if $\delta\mathcal{R} = 0$. According to the second Bianchi identity

$$\delta\mathcal{R} = -d(Ric(g));$$

we can express it in local coordinates $D_k \mathcal{R}_{ij} = D_i \mathcal{R}_{kj}$. In this case, $Ric(g)$ is a Codazzi tensor . In particular, any Riemannian manifold with Ricci parallel tensor $Dr = 0$ (i.e. $D_i \mathcal{R}_{kj} = 0$ ) has an harmonic curvature. Indeed, the Levi-Civita connexion $D$ is a Yang-Mills connexion on the tangent bundle of $M$ . In this way, the connexion $D$ is a critical point of the Yang-Mills functional

$$\mathcal{YM}(\nabla) = \frac{\infty}{\in} \int_{\mathcal{M}} ||\mathcal{R}^\nabla|| \lceil \sqsubseteq,$$

where $R^\nabla$ is the curvature associated to the connexion $\nabla$ . Notice that the Riemannian curvature must be harmonic for all Einstein manifolds and for all conformally flat manifolds with constant scalar curvature; this result can be deduced from the orthogonal decomposition of the curvature tensor. Moreover, the condition of the curvature harmonicity is, in a way, a generalisation of the Einstein condition of the metric ( $Ric(g)$ is a constant multiple of $g$ : $Ric(g) - \frac{R}{n}g = 0$ ). In particular, this fact shows that every Einstein metric must have a parallel Ricci tensor $Dr = 0$. In general, the latter property fails for a conformally flat metric with constant scalar curvature, as we shall see below. However, in the 3-dimensional case, we get an identity between harmonic curvature metrics and conformally flat metrics with constant scalar curvature.

Moreover, One asked about the existence of compact metric with harmonic curvature and non parallel Ricci curvature, as it was wellknown in the non compact case.

As response, Derdzinski has given examples of such metrics with $Dr \neq 0$ , [De]. The corresponding manifolds are bundles with fibres $N$ over the circle $S^1$ (parametrized by arc length t ) with allowed warped metrics $dt^2 + h^{4/n}(t)g_0$ on the product $S^1 \times N$. Here, $(N, g_0)$ is an Einstein manifold of dimension $n \geq 3$, and the function $h(t)$ on the prime factor is a periodic solution of the ODE, established by Derdzinski

$$h'' - \frac{nR}{4(n-1)}h^{1-4/n} = -\frac{n}{4}Ch \quad \text{for some constant} \quad C > 0. \quad (1)$$

This function must be non constant, otherwise the corresponding metric has a parallel Ricci tensor.

Notice that the manifolds $S^1 \times N$ are not conformally flat, unless $(N, g_0)$ has constant sectional curvature. In fact, under the hypothesis of non negative sectional curvature on a compact Riemannian manifold with harmonic curvature, the Ricci tensor is always parallel .

# 3 Curvature property of the pseudo-cylindric metrics

Let $(S^n, d\xi^2)$ be the standard sphere with radius 1 . Consider metrics of the conformal class $g \in [d\xi^2]$ on the domain $S^n - \Lambda_k$, where $\Lambda_k$ is a finite point-set of $S^n$.

We are interested on the complete metrics $g$ conformal to $d\xi^2$ with (posi-tive) constant scalar curvature $R(g)$.

For $k = 2$, there is a conformal diffeomorphism between $S^n - \{p_1, p_2\}$ and $(S^1 \times S^{n-1}, dt^2 + d\xi^2)$, where $S^1$ is the circle of length $T$. The Yamabe metrics on $(S^1 \times S^{n-1}, dt^2 + d\xi^2)$, are called pseudo-cylindric metrics. There are metrics of the form $g = u^{\frac{4}{n-1}}(dt^2 + d\xi^2)$ the $C^\infty$ function $u$ is a non constant positive solution of the Yamabe equation, [Ch].

Furthermore, Derdzinski established a classification of the compact n-dimensional Riemannian manifolds $(M_n, g)$, $n \geq 3$, with harmonic curvature. If the Ricci tensor $Ric(g)$ is not parallel and has less than three distinct eigenvalues at each point, then $(M, g)$ is covered isometrically by a manifold
$$(S^1(T) \times N, dt^2 + h^{4/n}(t)g_0),$$
where the non constant positive periodic solutions $h$ verify the equation (1). Here $(N, g_0)$ is a (n-1)- dimensional Einstein manifold with positive (constant) scalar curvature.
We get the following result, which gives another metric sharing this curvature property.

**Theorem 1** *Consider the product manifold $(S^1(T) \times S^{n-1}, g)$, $g = dt^2 + d\xi^2$ where $(S^1, dt^2)$ is the circle of length $T$ and $(S^{n-1}, d\xi^2)$ is the standard sphere with radius 1. Under the condition*
$$T > T_1 = \frac{2\pi}{\sqrt{n-2}},$$
*on the circle length, the Riemannian curvatures of the associated pseudo-cylindric metrics $g_c = u_c^{\frac{4}{n-2}} g$ are harmonic and their Ricci tensors are non parallel.*
*Moreover, any pseudo-cylindric metric may be identified to a Derdzinski metric up to a conformal transformation.*

## 3.1 Proof of Theorem 1

The local coordonates are $(x^0, x^1, ...., x^{n-1})$, and $x^0 = t$ is the coordinate corresponding to the circle factor $S^1$. We can write any metric tensor in the conformal class $[g_0]$ : $\bar{g} = u^{\frac{4}{n-2}}g,$ where the $C^\infty$ function $u$ is defined on the circle. Then we obtain the expressions in the local coordonates system

$$\bar{g}_{00} = u^{\frac{4}{n-2}}, \quad \bar{g}_{0i} = 0, \quad \bar{g}_{ij} = u^{\frac{4}{n-2}}(g_0)_{ij},$$

where $i,j = 1,2,...,n-1.$

We get the Christoffel symbols $\bar{\Gamma}^i_{jk}$ associated to the metric $\bar{g}$.

$$\bar{\Gamma}^0_{jk} = -\frac{2}{n-2}\frac{u'}{u}\bar{g}_{jk}, \quad \bar{\Gamma}^i_{j0} = \frac{2}{n-2}\frac{u'}{u}\delta^i_j, \quad \bar{\Gamma}^0_{00} = \frac{2}{n-2}\frac{u'}{u} \text{ and } \bar{\Gamma}^i_{00} = 0$$

The Ricci tensor components are

$$\bar{\mathcal{R}}_{ij} = \mathcal{R}_{ij} - \frac{2\nabla_{ij}u}{u} + \frac{2n}{n-2}\frac{\nabla_i u \nabla_j u}{u^2} - \frac{2}{n-2}\frac{|\nabla u|^2 + u\Delta u}{u^2}g_{ij}.$$

We have
$$\bar{\mathcal{R}}_{0i} = 0, \quad \bar{\mathcal{R}}_{00} = 2\frac{n-1}{n-2}[\frac{u''}{u} + \frac{u'^2}{u^2}].$$

The associated scalar curvatures are

$$R(g_0) = (n-1)(n-2) \quad and \quad \bar{R}(u^{\frac{4}{n-2}}g_0) = n(n-1).$$

Let $D$ denotes the riemannian connexion associated to the metric $g$. Since $u$ should be a non constant periodic function, a trivial calculus gives

$$D_0\bar{\mathcal{R}}_{00} = \frac{d\bar{\mathcal{R}}_{00}}{dt} \neq 0.$$

The Ricci curvature of the pseudo-cylindric metrics is not parallel, except for the cylindric one.

Return now to the Derdzinski metrics descrived above. We remark that, only in the locally conformally flat case of this product, the factor $N$ may be identified with the standard sphere $S^{n-1}$. Hence, these metrics must be conformally flat and will be conformal to the cylindric metric, as S.T. Yau [Y] has shown

**Lemma 1**
*Any warped metric $dt^2 + f^2(t)g_0$ on the product Riemannian manifold*

$(S^1(T) \times S^{n-1}, dt^2 + g_0)$ *must be conformally flat and conformal to a Riemannian metric product* $d\theta^2 + d\xi^2$ *where* $\theta$ *is a $S^1$-parametrisation with length* $\int_{S^1} \frac{dt}{f(t)}$ .

Thus, any Derdzinski metric may be identified with a pseudo-cylindric metric up to conformal diffeomorphism, let $F$ . Then the metrics are relied

$$dt^2 + f^2(t)d\xi^2 = F^*(u^j{}_c{}^{\frac{4}{n-2}}(dt^2 + d\xi^2)),$$

where $u^j{}_c$ are the Yamabe solutions belonging to the (same) conformal class.

Indeed, for the metric $dt^2 + f^2(t)d\xi^2$ we can write

$$dt^2 + f^2(t)d\xi^2 = f^2(t)[(\frac{dt}{f})^2 + d\xi^2].$$

After a change of variables and by using the conformal flatness of the product metric, we get

$$dt^2 + f^2(t)g_0 = \phi^2(\theta)[d\theta^2 + d\xi^2]$$

which is conformally flat.
To see that, it suffices to remark that any manifold carying a warped metric product $(S^1 \times N, dt^2 + h^{4/n}(t)g_0)$ with harmonic curvature is not conformally flat unless $(N, g_0)$ is a space of constant curvature. This manifold must be locally conformally equivalent to the trivial product $S^1 \times N$.
This product will be conformally flat only if $N$ has constant sectional curvature.

## 3.2 Remark concerning the parallelism of the Ricci tensor

Parallelism property of the Ricci tensor have an interest particularly for the conformally flatness case. Indeed, consider any Riemannian manifold $(M, g)$ with a parallel Ricci tensor. This implies in particular, that its Weyl tensor is harmonic: $\delta W = 0$ (in local coordonates $D_k W^h_{ijk} = 0$). For a metric $\tilde{g}$ in the conformal class and with harmonic Weyl tensor $[g]: \tilde{g} = e^{2\rho}g$, we get the following equality

$$\tilde{\delta}\tilde{W} = \delta W - \frac{1}{2}(n-3)W(\nabla\rho, ., ., .).$$

Since, the corresponding Ricci tensor is parallel, we then obtain in local coordonates
$$(n-3)W^l_{ijk}D_l\rho = 0.$$
Thus, when $(M,g)$ is not conformally flat, then $\rho = 0$ necessarily (see "The Ricci calculus" of Schouten).

# 4 Ricci tensor of the asymptotic pseudo-cylindric metrics

In this section we consider singular Yamabe metrics on the manifold $M = S^n - \Lambda_k$ where $k > 2$. These metric are complete and have a (positive) constant scalar curvature. Which property must it have their Ricci tensor ? We consider for that the estimate of the positive solutions $u(t,\xi)$, on the cylinder of the equation

$$4\frac{n-1}{n-2}\Delta_{g_0}u + R(g_0)u - R(g)u^{\frac{n+2}{n-2}} = 0,$$

$$g = u^{\frac{4}{n-2}}g_0 \quad complete \ on \quad S^n - \Lambda_k \quad and \quad R(g) = constant > 0 \ .$$

By using a reflection argument, it is known [M-P] that any solution $u(x)$ of the above equation, with a singularity at $p$, is asymptotic to one of the pseudo-cylindric functions near the point $p$. In fact, any such solution $u(x)$ corresponds to a solution $u(t,\xi)$ on a domain of $I\!R \times S^{n-1}$.
Moreover, the corresponding pseudo-cylindric metric is unique.

Let $\bar{g} = [u(t,\xi)]^{\frac{4}{n-2}}g_0$ be a Yamabe metric conformal to the standard metric and $u^{j,l}{}_c(t)$ be the pseudo-cylindric solution at a point $p_l \in \Lambda_k$. The index $j$ is corresponding to one solution having period $T$.
More precisely, we have the following .

**Lemma 2**
As $t \to \infty$ we get the following estimate of the Yamabe solution on $S^n - \Lambda_k$ near the point $p_l$

$$u(t,\xi) = u^{j,l}{}_c(t) + u^{j,l}{}_c(t)\mathcal{O}(e^{-\beta t}) \ ,$$

where $u^{j,l}{}_c(t)$ is the corresponding pseudo-cylindric solution and $\beta$ is a positive constant.

In the other words, this solution $u(t,\xi)$, asymptotically closed to a pseudo-cylindric function $u^j{}_c(t)$, around a singular point $p_l$, can be written as $t \to \infty$

$$u(t,\xi) = u^j{}_c(t) + e^{-\alpha t} v(t,\xi), \qquad (2)$$

where $\alpha$ is a positive constant, only depending on the cylindric bound and $v(t,\xi)$ is bounded.

So, the corresponding metrics on $S^n - \Lambda_k$ are complete and asymptotic at each singular point to a (unique) pseudo-cylindric metric.

Thus, we may estimate the Yamabe metrics on the manifold $S^n - \Lambda_k$, where $\Lambda_k$ is a finite point-set of the standard sphere $(S^n, g_0)$. A Yamabe metric $\bar{g}$ on $S^n - \Lambda_k$ conformal to the standard one then is asymptotic to a (unique) pseudo-cylindric metric $g$ around a singular point $p_l$.

We get the following result

**Theorem 2**
*Any complete metric $\bar{g}$ of positive constant scalar curvature which is conformal to the standard metric $g_0$ on the manifold $S^n - \Lambda_k$, where $\Lambda_k = \{p_1, p_2, ..., p_k\}$ $k > 2$, has a harmonic curvature and non parallel Ricci tensor, except for the standard one.*

## 4.1 Proof of theorem 2

We use the asymptotic property of these metrics (Lemma 2). In particular, we get

$$u(t,\xi) = u^j{}_c(t) + e^{-\alpha t} v(t,\xi).$$

Notice that, this manifold $(S^n - \Lambda_k, \bar{g})$ is locally conformally flat. The condition of constant scalar curvature implies that its Riemannian curvature is harmonic. We get the converse only for the dimension $n = 3$.

We proceed as in the proof of Theorem 1. Consider a local chart $(t, \xi)$ of a point in $S^n - \Lambda_k$, here $t = x^0$ and $\xi = (x^1, x^2, ..., x^{n-1})$.
In this local
coordonates system, we obtain the components of the tensor metric

$$\bar{g}_{00} = u^{\frac{4}{n-2}} \quad , \quad \bar{g}_{0i} = 0 \quad , \quad \bar{g}_{ij} = u^{\frac{4}{n-2}}(g_0)_{ij}.$$

Then we get components of the Ricci tensor

$$\bar{\mathcal{R}}_{0i} = -\frac{2\nabla_{i0}u}{u} + \frac{2n}{n-2}\frac{1}{u^2}\frac{\partial u}{\partial x_i}\frac{\partial u}{\partial t},$$

$$\bar{\mathcal{R}}_{00} = 2\frac{n-1}{n-2}[\frac{1}{u}\frac{\partial^2 u}{\partial t^2} + \frac{1}{u^2}(\frac{\partial u}{\partial t})^2].$$